\documentclass[12pt]{amsart}
\usepackage{amsmath,amssymb, xypic,verbatim,amscd,color}
\usepackage{graphicx}
\usepackage{colordvi}
\usepackage{times}

\newtheorem{theorem}{Theorem}[section]

\newtheorem{lemma}[theorem]{Lemma}

\newtheorem{definition}[theorem]{Definition}

\theoremstyle{definition}

  \newcommand{\sig}{\sigma}

  \newcommand{\gam}{\gamma}
  \newcommand{\lam}{\lambda}

  \newcommand{\Ker}{\operatorname{Ker}}

 \newcommand{\Exp}{\operatorname{Exp}}

\newcommand{\beqn}{\begin{equation}}
\newcommand{\eeqn}{\end{equation}}

 \newcommand{\fb}{\mathfrak{b}}
 \newcommand{\fg}{\mathfrak{g}}
\newcommand{\fp}{\mathfrak{p}}
 
\newcommand{\fl}{\mathfrak{l}}

 \newcommand{\fs}{\mathfrak{s}}
 \newcommand{\ft}{\mathfrak{t}}

\newcommand{\bc}{\mathbb{C}}

 \newcommand{\cl}{\mathcal{L}}

\begin{document}
\title[equivariant cohomology of some Springer fibers]
{An algebro-geometric realization
of equivariant cohomology of some Springer fibers}

\author{Shrawan Kumar}
\author{Claudio Procesi}
\footnote{S.K. was supported by NSF  grants.}

 \maketitle

\begin{abstract} We give an explicit affine algebraic variety whose
coordinate ring is isomorphic (as a $W$-algebra) with the equivariant cohomology
of some Springer fibers. \end{abstract}

\section{Introduction}

Let $G$ be a
connected simply-connected semisimple  complex algebraic group with a Borel subgroup
$B$ and a maximal torus $T\subset B$. Let $P\supseteq B$ be a (standard)
 parabolic subgroup of $G$. Let $L\supset T$ be the Levi subgroup of $P$ and let $S$
 be the connected center of $L$ (i.e., $S$ is the identity component of the center of $L$).
 Then, $S\subset T$. We denote
 the Lie algebras of $G,T,B,P,L,S$  by the corresponding Gothic characters: $\fg,\ft,
 \fb,\fp,\fl,\fs$ respectively.
 Let $W$ be the Weyl group of $G$ and $W_L \subset W$ the Weyl group of $L$. Let
 $\sigma = \sigma_\fl$ be
 a principal nilpotent element of $\fl$. Let $X=G/B$ be the full flag variety of $G$
 and let $X_\sigma\subset X$ the Springer fiber corresponding to the nilpotent element
 $\sigma$ (i.e., $X_\sigma$ is the
 subvariety of $X$
 fixed under the left multplication by $\Exp \sigma$ endowed with the reduced
 subscheme structure). Observe that $S$ keeps
 the variety $X_\sigma$ stable under the left multiplication of $S$ on $X$.

 \begin{definition} \label{def} Let $Z_\fl$ be the reduced
 closed subvariety of $\ft\times \ft$ defined by:
 \[Z_\fl :=\{(x,wx): w\in W, x\in \fs\}.\]
 Since $Z_\fl$ is a cone inside $\ft\times \ft$, the affine coordinate ring
 $\Bbb C[Z_\fl]$ is a non-negatively graded algebra. Moreover, the projection
 $\pi_1: Z_\fl \to \fs$ on the first factor gives rise to a $S(\fs^*)$-algebra
 structure on $\Bbb C[Z_\fl]$.
Also, define an action of $W$ on $Z_\fl$ by:
\[v\cdot (x,wx)= (x,vwx), \,\,\text{for}\, x\in \fs, v,w\in W.\]
This action gives rise to a $W$-action on $\Bbb C[Z_\fl]$, commuting with the
$S(\fs^*)$ action
 on $\Bbb C[Z_\fl]$.

 In fact, even though we do not need it, $W$ is precisely the automorphism group
 of 
 $\Bbb C[Z_\fl]$ as $S(\fs^*)$-algebra.

 For $\fp = \fb$, the Levi subalgebra $\fl$ coincides with $\ft$, $\sigma_\ft=0$ and
 $X_\sigma = X$. In this case, $\fs= \ft$ and we abbreviate $Z_\ft$ by $Z$. Clearly, $Z_\fl$
 (for any Levi subalgebra $\fl$) is a
 closed subvariety of $Z$.
\end{definition}

The following theorem  is our main result.

 \begin{theorem}  \label{main} With the notation as above, assume that the canonical restriction
 map $H^*(X)\to H^*(X_\sigma)$ is surjective, where $H^*$ denotes the singular
 cohomology with complex coefficients. Then, there is a graded
 $S(\fs^*)$-algebra isomorphism
 \[\phi_\fl: \Bbb C[Z_\fl] \to H^*_S(X_\sigma),\]
 where $H^*_S$ denotes the $S$-equivariant cohomology with complex coefficients.

 Moreover, the following diagram is commutative:
 \beqn \label{e1} \xymatrix{ \Bbb C[Z]\ar[d] \ar[r]^{\phi_\ft} & H_T^*(X)\ar[d]\\
   \Bbb C[Z_\fl] \ar[r]^{\phi_\fl} & H_S^*(X_\sigma),  } \eeqn
   where the vertical maps are the canonical restriction maps.

 In particular, we get an isomorphism of graded algebras
 \[\phi_\fl^o: \Bbb C\otimes_{S(\fs^*)}\, \Bbb C[Z_\fl] \to H^*(X_\sigma ),\]
 making the following diagram commutative:
 \beqn \label{e2}\xymatrix{ \Bbb C \otimes_{S(\ft^*)}  \Bbb C[Z]\ar[d]\ar[r]^{\phi_\ft^o} & H^*(X)\ar[d]\\
  \Bbb C \otimes_{S(\fs^*)} \Bbb C[Z_\fl] \ar[r]^{\phi_\fl^o} & H^*(X_\sigma),  } \eeqn
  where the vertical maps are the canonical restriction maps and $\Bbb C$ is
  considered as  a $S(\fs^*)$-module under the evaluation at $0$.

 Moreover, the isomorphism $\phi_\fl^o$ is $W$-equivariant under the Springer's
 $W$-action on $H^*(X_\sigma )$ and the $W$-action on
 $ \Bbb C\otimes_{S(\fs^*)}\, \Bbb C[Z_\fl]$ induced from the $W$-action on
 $\Bbb C[Z_\fl]$ defined above.
 \end{theorem}

\vskip2ex

\noindent
{\bf Acknowledgements.} We thank Eric Sommers for  the example given in
Remark  2.3  (3).

\section{Proof of the Theorem}

Before we come to the proof of the theorem, we need the following lemma. (See,
e.g., [C, Theorem 2].)

 \begin{lemma}\label{lemma1}
For any $w\in W$, there exists a unique $w'\in W_L$
such that
  \[
w'w B \in X^S_\sigma \subset X.
  \]

  Moreover, this induces a bijection
  \[
W_L \backslash W \leftrightarrow X^S_\sigma .
  \]
\end{lemma}

We also need the following simple (and well known) result.

\begin{lemma} \label{lemma2} Let $S=S(V^*)$ be the symmetric algebra for a finite dimensional
vector space $V$ and let $M,N,R$  be three $S$-modules. Assume that $N$ and $R$
 are $S$-free of the same finite rank and $M$ is a $S$-submodule of $R$. Then,
 any surjective $S$-module morphism $\phi:M \to N$ is an isomorphism.
 \end{lemma}

We now come to the proof of the theorem.

\vskip1ex
\noindent
{\bf Proof of the theorem.}
Consider the equivariant Borel homomorphism
  \[
\beta : S(\ft^*) \to H_T(X)
  \]
obtained by $\lam\mapsto c_1(\cl_\lam )$, where $\lam\in\ft^*$ and
$c_1(\cl_\lam )$ is the $T$-equivariant first Chern class of the line
bundle $\cl (\lam )$ on $X$ corresponding to the character $e^\lam$, and extended
as a graded algebra homomorphism.
This gives rise to an algebra homomorphism
  $$
\chi : \bc [\ft\oplus\ft ] \simeq S(\ft^*) \otimes S(\ft^*) \to H_T(X),\,\,\,
p\otimes q \mapsto p\cdot \beta (q),
  $$
where $p\cdot$ denotes the multiplication in the $T$-equivariant cohomology
by $p\in S(\ft^*)\simeq H_T(pt)$.  It is well known that $\chi$ is
surjective.  Moreover, both the restriction maps
  \[
H_T(X) \twoheadrightarrow H_S(X) \twoheadrightarrow H_S(X_\sig )
  \]
are surjective; this follows since both the spaces $X$ and $X_\sig$
have cohomologies  concentrated in even degrees (cf. [DLP]).  (Use the degenerate
Leray-Serre
spectral sequence and the assumption that the restriction map
$H^*(X) \to H^*(X_\sigma)$ is surjective.)

Consider the canonical surjective map $\theta : \bc [\ft\oplus\ft ]
\twoheadrightarrow \bc [Z_\fl ]$.  Then, of course,
  \beqn \label{e3}
\Ker \theta = \biggl\{ \sum_i p_i\otimes q_i : p_i,q_i \in S(\ft^*)\,\text{and}\,
\sum_i p_i(x)\, q_i
(wx) =0, \text{ for all $x\in\fs$ and $w\in W$}\biggr\}.\eeqn
  We claim that
\beqn \label{e4} \Ker \theta\subset \Ker \gam,\eeqn
where $\gam$ is the
composite  map
  \[
\bc [\ft\oplus\ft ] \overset{\chi}\to H_T(X) \to H_S(X_\sig ).
  \]
Since $X_\sigma$ has cohomologies only in even degrees, by the degenerate
Leray-Serre spectral sequnce, $H_S(X_\sig )$ is a
free $S(\fs^*)$-module. In particular,    by the Borel-Atiyah-Segal Localization
Theorem (cf. [AP, Theorem 3.2.6]),
  \[
H_S(X_\sig ) \hookrightarrow H_S(X^S_\sigma ).
  \]
Thus, to prove the claim ~ (\ref{e4}),  it suffices to prove that for any $\sum_i
p_i\otimes q_i\in \Ker\theta$,
  \[
\gam \biggl( \sum_i p_i\otimes q_i\biggr)_{|_{X^S_\sigma}} \equiv 0.
  \]
It is easy to see that the Borel homomorphism $\beta$ restricted to the
$T$-fixed points $X^T$ satisfies:
  \[
\beta (q) (wB) = wq,\; \text{ for any } q\in S(\ft^*) \,\text{and}\, w\in W.
  \]
Thus, for any $w\in W$,
  \[
\gam \biggl( \sum_i p_i\otimes q_i\biggr) (w'wB) = \biggl(
\sum_i(p_i)(w'wq_i)\biggr)_{|_\fs},
  \]
  where $w'$ is as in Lemma ~\ref{lemma1}.
From the description of $\Ker \theta$ given in ~(\ref{e3}), we thus get that the
claim ~(\ref{e4})  is true.
Hence, the map $\theta$ descends to a surjective $S(\fs^*)$-module map
  \[
\phi_\fl : \bc [Z_\fl ] \twoheadrightarrow H_S(X_\sig ).
  \]

Again using the Localization Theorem, the free $
S(\fs^*)$-module $H_S(X_\sig )$ is   of rank $=\#\, W_L \backslash W$, since $\#\ X^S_\sigma=
\#\, W_L \backslash W$ by Lemma ~\ref{lemma1}.  Also, the projection
on the first factor $\pi_1: Z_\fl \to \fs$ is a finite morphism with all its
fibers of cardinality $\leq \#\, W_L \backslash W$.  To see this, consider the
 surjective morphism $\alpha:\fs \times W/W_L \to Z_\fl, \, (x, wW_L)\mapsto
 (x,wx)$. Then, $\pi_1\circ \alpha : \fs \times W/W_L \to \fs$ is again
 the projection on the first factor, which is clearly a finite morphism and hence so
 is $\pi_1$.

 Now, taking $M=\Bbb C[Z_\fl], N=H_S^*(X_\sigma), R=\Bbb C[\fs \times W/W_L]$
 and $V=\fs$ in Lemma ~\ref{lemma2}, we get that $\phi_\fl$ is an isomorphism,
 where the inclusion $M \subset R$ is induced from the surjective morphism
 $\alpha:\fs \times W/W_L \to Z_\fl$.

 The commutativity of the diagram ~(\ref{e1})  clearly follows from the above proof.

 Since $H^*(X_\sigma)$ is concentrated in even degrees, by the degenerate
 Leray-Serre spectral sequence, we get that
 \[H^*(X_\sigma) \simeq \Bbb C \otimes_{S(\fs^*)} H_S^*(X_\sigma).\]
 From this the `In particular' part of the theorem follows.

 From the definition of the map $\phi_\ft$, it is clear that $\phi_\ft^o$ is
 $W$-equivariant with respect to the action of $W$ on  $\Bbb C \otimes_{S(\ft^*)}
 \Bbb C[Z]$ induced from the action of $W$ on $\Bbb C[Z]$ as defined in Definition
 ~\ref{def} and the standard action of $W$ on $H^*(X)$. Moreover, the restriction map
 $H^*(X)\to H^*(X_\sigma)$ is $W$-equivariant with respect to the Springer's
 $W$ action on $H^*(X_\sigma)$ (cf. [HS, \S2]). Thus, the $W$-equivariance of
 $\phi_\fl^o$ follows from the commutativity of the diagram ~(\ref{e2}). This
 completes the proof of the theorem. \qed

\vskip2ex

\noindent
  {\bf Remark 2.3.} \label{remark}   (1) By the Jordan block decomposition, any nilpotent element
  $\sigma \in sl(N)$ (up to conjugacy)
  is a regular nilpotent element in a standard Levi subalgeabra $\fl$ of $sl(N)$.
  Moreover,
the canonical restriction
 map $H^*(X)\to H^*(X_\sigma)$ is surjective in this case. In fact, as proved by
 Spaltenstein [S], in this case there is a paving of $X$ by affine spaces as
 cells such that $X_\sigma$ is a closed union of cells (cf. also [DLP]).
 Thus, the above theorem, in particular, applies to any nilpotent element
  $\sigma$ in
  any special linear Lie algebra $sl(N)$.

  (2) A certain variant (though a less precise version) of our Theorem ~\ref{main}
  for $\fg=sl(N)$ is obtained by Goresky-MacPherson [GM, Theorem 7.2].

  (3) For a general semisimple Lie algebra $\fg$,  it is not true that the restriction map
  $H^*(X)\to H^*(X_\sigma)$ is surjective  for any regular nilpotent element in a
  Levi subalgebra $\fl$. Take, e.g., $\fg$ of type $C_3$ and $\sigma$ corresponding
  to the Jordan blocks of size $(3,3)$. In this case, the centralizer of $\sigma$ in the
  symplectic group Sp$(6)$ is connected and $X_\sigma$ is two dimensional. The cohomology of
  $X_\sigma$ as a $W$-module is given as follows:

  Of course, $H^0(X_\sigma)$ is the one dimensional  trivial  $W$-module;
  $H^2(X_\sigma)$ is the sum of the three dimensional reflection representation with a
  one dimensional representation; and $H^4(X_\sigma)$ is a three dimensional
  irreducible representation.

\bibliographystyle{plain}
\def\noopsort#1{}

\vskip6ex

\noindent
Addresses:

 S.K: Department of Mathematics, University of North Carolina, Chapel Hill,
 NC 27599-3250, USA\\
\noindent(email: shrawan@email.unc.edu)

C.P.: Dipartimento di Matematica, Sapienza Universit� di Roma, Piazzale Aldo Moro 5,
 00185 Roma, Italy.
\\
\noindent(email:procesi@mat.uniroma1.it)

\end{document}